\documentclass[reqno,12pt]{amsart}%,draft
\usepackage[mathscr]{eucal}
\usepackage{amsmath}
\usepackage{amssymb}
\usepackage{latexsym}
\usepackage{amsthm}
\usepackage{color}
\newtheorem{thm}{\indent Theorem}[section]

\newtheorem{lem}[thm]{\indent Lemma}
\newtheorem{prop}[thm]{\indent Proposition}
\newtheorem{dfn}{{\indent\bf Definition}}[section]
\newtheorem{rmk}{{\indent\bf Remark}}[section]

\newcommand{\lmx}{\left(\begin{matrix}}
\newcommand{\rmx}{\end{matrix}\right)}
\newcommand{\ldt}{\left|\begin{matrix}}
\newcommand{\rdt}{\end{matrix}\right|}

\newcommand{\td}{\tilde}

\newcommand{\nnm}{\nonumber}
\newcommand{\bbr}{{\mathbb R}}

\newcommand{\be}{\begin{equation}}
\newcommand{\ee}{\end{equation}}

\makeatletter
\numberwithin{equation}{section}
\begin{document}
\title [affine maximal type hypersurfaces]{{A class of new complete affine maximal type hypersurfaces $^\star$}}
\author [Sun and Xu] {Yalin Sun \quad\quad Ruiwei Xu$^*$ }
\address{School of Mathematics and Statistics,
\newline \indent Henan Normal University, Xinxiang 453007, P. R. China
\newline \indent ylsun2024@126.com;\; rwxu@htu.edu.cn}
\date{}
\footnotetext{$\star$\, This project is supported by NSFC of China, Grant Number
 11871197 and Key projects of the Natural Science Foundation of Henan Province (No. 252300421303). \\ \indent$*$ Corresponding author.}
\maketitle
\renewcommand{\baselinestretch}{1}
\renewcommand{\arraystretch}{1.2}
\catcode`@=11 \@addtoreset{equation}{section} \catcode`@=12
\makeatletter
\renewcommand{\theequation}{\thesection.\arabic{equation}}
\@addtoreset{equation}{section} \makeatother
%%%%%%%%%%%%%%%%%%%%%%%%%%%%%%%%%%%%%%%%%%%%%%%
\noindent {\bf Abstract:} In this paper we classify a kind of special Calabi hypersurfaces with negative constant sectional curvature in Calabi affine geometry. Meanwhile, we find a class of new Euclidean complete and Calabi complete affine hypersurfaces, which satisfy the affine maximal type equation and the Abreu equation with
negative constant scalar curvatures.

%%%%%%%%%%%%%%%%%%%%%%%%%%%%%%%%%%%%%%%%%%%%%%
\vskip 2pt\noindent {\bf
2000 AMS Classification:} Primary 53A15; Secondary 35J60, 53C24.
%%%%%%%%%%%%%%%%%%%%%%%%%%%%%%%%%%%%%%%%%%%%
 \vskip 0.1pt

\section{Introduction}

In affine differential geometry one often encounters fourth order
non-linear PDEs which are far from being well understood.
Consequently, new and significant efforts are required for their
investigation. In this
paper we find a class of new complete solutions to the following nonlinear,
fourth order PDE for a
convex function $f(x)$ on a convex domain $\Omega$ in
$\mathbb{R}^{n}.$ The equation can be written as
\begin{equation}\label{1.1}\sum_{i,j=1}^nF^{ij} w_{ij} =-L^\sharp,
\;\;\;w =\left[\det\left(\tfrac{\partial^{2}f}{\partial
x_{i}\partial
x_{j}}\right)\right]^{a},\end{equation}
where $L^\sharp:\Omega\rightarrow \mathbb{R}$ is some given function, $(F^{ij})$ denotes the cofactor matrix of the Hessian matrix
$\left(\tfrac{\partial^{2}f}{\partial x_{i}\partial x_{j}}\right)
$ and $a\neq 0$ is a constant. In the case when $a = -\tfrac{n+1}{n+2}$, the PDE (\ref{1.1}) is the equation for  equiaffine hypersurfaces with affine mean curvature $\frac{L^\sharp}{n(n+1)}$ (see p.86 in \cite{LSZH}).
In particular, when $L^\sharp=0$, it is
the equation for \emph{affine maximal hypersurfaces}. In the case when $a=-1$, the PDE (\ref{1.1}) is called the
\emph{Abreu equation}, which appears in the study of the differential geometry of toric varieties (\cite{A, D-2}), where $L^\sharp$ is the scalar curvature of the K\"{a}hler metric.
As in \cite{D}, we call the PDE (\ref{1.1}) with $L^\sharp=0$, i.e.,
\begin{equation}\label{1.2}\sum_{i,j=1}^nf^{ij} w_{ij} =0,
\;\;\;w =\left[\det\left(\tfrac{\partial^{2}f}{\partial
x_{i}\partial
x_{j}}\right)\right]^{a},\end{equation}
as \emph{affine maximal type equation} in this paper.

Quadratic polynomial is the trivial solution to affine maximal type equation (\ref{1.2}). Bernstein problem for affine maximal type equation has been a core and important problem in affine geometry since the 1970s.
In 2000, Trudinger-Wang \cite{TW} first solved Chern's conjecture on affine maximal surfaces, and they proved that
the entire solution of the PDE (\ref{1.2}) on $\mathbb R^2$ with $a>0$ must be
a quadratic polynomial in \cite{TW1}. Later, Li-Jia \cite{LJ} proved that the entire solution of the PDE (\ref{1.2}) on $\mathbb R^2$ with $a\leq -\frac{3}{4}$ must be
a quadratic polynomial, which gave a new analytic proof
for Chern's conjecture. For higher dimensions, Jia and
Li \cite{JL} proved the uniqueness
of the entire solutions to the PDE (\ref{1.2}) when $|a|$ is large enough, see also \cite{XLL}. Up to now the higher dimensional Chern's conjecture is still open.

Besides Euclidean completeness, one can consider another possibility of assuming the completeness with respect to the Calabi metric. In this aspect, the first result that was obtained by Li and Jia in \cite{L-J-4}:  if the affine maximal hypersurface is complete with respect to the Calabi metric, then, for dimension $n=2$ or $n=3$,  the hypersurface must be an elliptic paraboloid. If the affine maximal type hypersurfaces are complete with respect to the Calabi metric, there is the following more general result, see Theorem 5.8.1 in \cite{LXSJ} or (2) in Theorem 1.1 of \cite{D2}. The other Bernstein type results of the PDE (\ref{1.2}), please see the papers \cite{CSL, LXSJ1,L,Xu-Cao}, etc.
\begin{thm} \cite{LXSJ}
For $n \geq 2$, if
\begin{equation}\label{1.3}
a \notin\left[ -\tfrac{1}{2}-\tfrac{n-1}{4\sqrt{n}}, -\tfrac{1}{2}+\tfrac{n-1}{4 \sqrt{n}} \right]
\end{equation}
and the Calabi metric is complete, then the solution to the PDE (\ref{1.2}) must be a quadratic function.
\end{thm}

On the opposite orientation, one hope to find counter-examples to Bernstein problem of the PDE (\ref{1.2}).
Recently, Du \cite{D} constructed various counter-examples explicitly and proved the following results about Euclidean complete affine maximal type hypersurfaces.
\begin{thm}\cite{D}
Assuming that $n\geq2$ and $a\in[-\tfrac{n-1}{n},0)$, there exist Euclidean complete affine maximal type hypersurfaces which are not elliptic paraboloid.
\end{thm}

In this paper we classify a kind of special Calabi hypersurfaces with negative constant sectional curvature in Calabi affine geometry. Meanwhile, we find new non-quadratic Euclidean complete and Calabi complete solutions of the PDE (\ref{1.2}) with $a=-\tfrac{n}{n+1}$ and $-\tfrac{1}{n+1}$.

\begin{thm}\label{thm1.3} Supposing that $n\geq2$ and $a=-\tfrac{n}{n+1}$ or $a=-\tfrac{1}{n+1}$, there exist non-quadratic Euclidean complete and Calabi complete affine maximal type hypersurfaces which are given by respectively
\be \label{1.4} f(x)=-\tfrac{1}{4}\ln{(x_1-\sum^n_{k=2}\tfrac{x^2_k}{2})},\;\;\;\;\;  x_1>\sum^n_{k=2}\tfrac{x^2_k}{2}, \ee
and
\be \label{1.5}
f(x)=-\tfrac{1}{4}\ln{x_1}+\sum^n_{k=2}\tfrac{x^2_k}{2x_1},\;\;\;\;\;  x_1>0.
\ee
They also satisfy the Abreu equation with scalar curvature $-4(n+1)$ and $-4n(n+1)$, respectively.
\end{thm}

Note that $a=-\tfrac{n}{n+1}$ is not in the scope of Theorem 1.2 of \cite{D}.
In case when $n=2$, we feel that examples (\ref{1.4}) and (\ref{1.5}) may be optimal non-quadratic Calabi complete solutions to the PDE (\ref{1.2}). Thus we pose the following conjecture.

\textbf{Conjecture.} Let $f$ be a smooth strictly convex function defined on a convex domain $\Omega$ in $\mathbb{R}^2$. If $f$ satisfies the PDE (\ref{1.2}) with $a\notin [-\tfrac{2}{3},-\tfrac{1}{3}]$, and the Calabi metric is complete, then $f$ must be a quadratic polynomial.

The paper is organized as follows. In section 2, we review relevant
materials for Calabi hypersurfaces and warped product. In section 3, we study a special Ejiri frame in Calabi affine geometry and its related characterization. In section 4, we give a complete classification of those Calabi hypersurfaces with such frame, see Theorem 4.1. Finally, in section 5, Theorem 1.3 is proved.

%%%%%%%%%%%%%%%%%%%%%%%
%%%%%%%%%%%%%%%%%%%%%%%%
\section{Preliminaries}

In this section, we shall recall some basic facts for the Calabi geometry, see \cite{Ca} or section 3.3.4 in \cite{LXSJ}. Let $f$ be a strictly convex $C^\infty$-function on a domain $\Omega \subset \mathbb{R}^n$.
Consider the graph hypersurface
\be\label{2.1} M^n:=\{(x_1,\cdots, x_{n+1})\;|\;x_{n+1}=f(x_1,\cdots,x_n), (x_1,\cdots,x_n)\in\Omega\}.\ee
This Calabi geometry can also be essentially realised as a very special relative affine geometry of the graph hypersurface (\ref{2.1})
 by choosing the so-called {\em Calabi affine normalization} with $$Y =(0,\cdots,0,1)\in \bbr^{n+1}$$
being the fixed relative affine normal vector field which we call the {\em Calabi affine normal}.

For the position vector $x=(x_1,\cdots,x_n,f(x_1,\cdots,x_n))$ we have the decomposition
\be \tfrac{\partial^2 x}{\partial x_i\partial x_j}=c_{ij}^k
x_*\tfrac{\partial }{\partial x_k}+f_{ij}Y,\ee
with respect to the bundle decomposition $\bbr^{n+1}=TM^n\oplus\bbr\cdot Y$, where the induced affine connection $c^k_{ij}\equiv 0$. It follows that the relative affine metric is nothing but the \emph{Calabi metric}
\be G=\sum \tfrac{\partial^2f}{\partial x_i\partial x_j}dx_idx_j.\ee
The Levi-Civita connection with respect to the metric $G$ has the Christoffel symbols
$$\hat \Gamma_{ij}^k=\tfrac{1}{2}\sum f^{kl}f_{ijl},$$
where and hereafter
$$f_{ijk}=\tfrac{\partial^3f}{\partial x_i\partial x_j\partial x_k}, \quad (f^{kl})=(f_{ij})^{-1}. $$
Then we can rewrite the Gauss structure equation as follows:\be x_{,ij}=\sum A_{ij}^kx_*\tfrac{\partial }{\partial x_k}+G_{ij}Y,\ee
where the subscript $ _{,ij} $ denotes the covariant derivative with respect to
the Levi-Civita connection $\hat\nabla$ of the Calabi metric. The Fubini-Pick tensor $A_{ij}^k$ and the Weingarten tensor satisfy
\be A_{ijk}=A_{ij}^lG_{kl}=-\tfrac{1}{2}f_{ijk},\quad B_{ij}=0,\ee
which means $A_{ijk}$ are symmetric in all indexes.
Classically, the tangent vector field
\be\label{2.5} T:=\tfrac{1}{n}\sum G^{kl}G^{ij}A_{ijk}\tfrac{\partial}{\partial x_l}=-\tfrac{1}{2n}\sum f^{kl}f^{ij}f_{ijk}\tfrac{\partial}{\partial x_l}\ee
is called the  \emph{Tchebychev vector field}, and the invariant function
\be\label{2.6} J:=\tfrac{1}{n(n-1)}\sum G^{il}G^{jp}G^{kq}A_{ijk}A_{lpq}=\tfrac{1}{4n(n-1)}\sum f^{il}f^{jp}f^{kq}f_{ijk}f_{lpq}\ee
is named as the  \emph{relative Pick invariant} of $M^n$.
The \emph{Gauss integrability conditions}  and the  \emph{Codazzi equations}  read
\be \label{2.7} R_{ijkl}=\sum f^{mh}(A_{jkm}A_{hil}-A_{ikm}A_{hjl}),\ee
\be \label{2.8} A_{ijk,l}=A_{ijl,k}.\ee
From (\ref{2.7}) we get the  \emph{Ricci tensor}
\be\label{2.9} R_{ik}=\sum f^{jl}f^{mh}(A_{jkm}A_{hil}-A_{ikm}A_{hjl}).\ee
Thus, the scalar curvature is given by
\be\label{2.10} R=n(n-1)J-n^2|T|^2.\ee

Let $A(n+1)$ be the group of $(n+1)$-dimensional affine transformations on $\bbr^{n+1}$. Then $A(n +1)=GL(n+1)\ltimes \bbr^{n+1}$, the semi-direct product of the general linear group $GL(n +1)$ and the group $\bbr^{n+1}$ of all the translations on $\bbr^{n+1}$. Define
\be\label{sa}
SA(n+1)=\{\phi=(M,b)\in A(n+1)=GL(n+1)\ltimes \bbr^{n+1};\ M(Y)=Y\}
\ee
where $Y=(0,\cdots,0,1)^t$ is the Calabi affine normal. Then the subgroup $SA(n+1)$ consists of all the transformations $\phi$ of the following type:
\begin{align}\label{sa1}
X:&=(X^j,X^{n+1})^t\equiv(X^1,\cdots, X^n,X^{n+1})^t\nnm\\
&\mapsto \phi(X):=\lmx a^i_j&0\\a^{n+1}_j&1\rmx X+b,\quad \forall\, X\in\bbr^{n+1},
\end{align}
for some $(a^i_j)\in A(n)$, constants $a^{n+1}_j$ ($j=1,\cdots,n$) and some constant vector $b\in \bbr^{n+1}$.

\begin{dfn}\cite{XL}\label{dfn1}
Two graph hypersurfaces $(x^i,f(x^i))$ and $(\td x^i,\td f(\td x^i))$, defined respectively in domains $\Omega,\td\Omega\subset\bbr^n$, are called Calabi-affine equivalence if they differ only by an affine transformation $\phi\in SA(n+1)$.
\end{dfn}

Next, we recall several concepts, see \cite{N1}. A subbundle $E\subset TM^n$ is called \emph{auto-parallel} if $\hat\nabla_XY\in E$ holds for all $X,\ Y\in E$. A subbundle $E$ is called \emph{totally umbilical} if there exists a vector field $H\in E^{\perp}$ such that
$$G(\hat\nabla_XY,Z)=G(X,Y)G(H,Z),$$
for all $X,\ Y\in E$ and $Z\in E^{\perp}$, here we call $H$ the mean curvature vector of $E$. If, moreover, $G(\hat\nabla_XH,Z)=0$ holds, we say that $E$ is \emph{spherical}. The following theorem of Hiepko \cite{N1} allows us to recognize when a manifold can be decomposed as a warped product manifold.
\begin{thm}
If the tangent bundle of the Riemannian manifold $M$ splits into an orthogonal sum $TM=E_0\bigoplus E_0^{\bot}$ of non-trivial subbundles such that $E_0$ is auto-parallel and $E_0^{\bot}$ is spherical, then the Riemannian manifold $M$ is locally isometric to a warped product $M_0\times_{\rho}M_1$ where $M_0$ and $M_1$ are integral manifolds of $E_0$ and $E_0^{\bot}$, respectively.
\end{thm}

%%%%%%%%%%%%%%%%%%%%%55
%%%%%%%%%%%%%%%%%%%%%%%%
\section{A special Ejiri frame and its related characterization}

Here we only consider the non-trivial solutions to the PDE (\ref{1.2}). Thus we assume that the Tchebychev vector field $T\neq 0$ on an open set $U\subset M^n$. Since our classification result Theorem 4.1 below is local in nature, in the following we restrict to consider Calabi hypersurface $M^{n}$ such that $U=M^{n}$.

Now, we fix a point $p \in M^n$. For subsequent purposes, we will review the well known construction of a typical orthonormal basis for $T_p M^n$, which was introduced by Ejiri and has been widely applied, and proved to be very useful for various situations, see e.g., \cite{CHM, HLV, LX,LV}.
Denote the set of unit vectors in $T_{p}M^{n}$ by
$$U_{p}M^{n}=\{u\in T_{p}M^{n}\mid G(u,u)=1\}.$$
Define a continuous function $F(u):=G(A_u u,u)$ on $U_{p}M^{n}$. Since $U_{p}M^{n}$ is compact and the Fubini-Pick tensor $A(p) \neq 0$, the function $F$ attains an absolute maximum at some direction $e_1 \in U_p M^n$ such that $\lambda:=F\left(e_1\right)>0$. Then the following lemma holds, see Lemma 3.1 in \cite{HLSV}.

\begin{lem}\label{lemma3.1}\cite{HLSV}\  There exists an orthonormal basis $\{e_1,\cdots,e_n\}$ of $T_pM^n$ such that the following hold:
\begin{equation}\label{3.1}
A_{e_1}e_{\alpha}=\lambda_{\alpha} e_{\alpha}, \quad 1\leq \alpha\leq n, \\
\end{equation}
where $\lambda_1>0,\; \lambda_1\geq 2\lambda_{i}$ for $2\leq i\leq n$. If, for some $i$,  $\lambda_i=\frac{1}{2}\lambda_{1}$, then $G(A_{e_i}e_i,e_i)=0. $
\end{lem}

Here we consider a special case of Ejiri frame: if, for any $2\leq i\leq n$, $\lambda_i=\frac{1}{2}\lambda_{1}$ on any point $p\in M$. To state the following proposition 3.2, we need the notion of the function $\theta$ introduced by Zhao-Cheng-Hu \cite{ZCH}. The function $\theta$ defined by
\begin{equation}\label{3.2}
\theta(p):=\max_{u\in U_pM^n}G(A_u u,u),\quad p\in M^n,
\end{equation}
is a continuous function on $M^n$, see Lemma 3.2 in \cite{ZCH}. Then we have the following characterization proposition about this kind of Ejiri frame.

\begin{prop}\label{prop3.2}\
 Let $x: M^n \rightarrow \mathbb{R}^{n+1}$ be a Calabi hypersurface. Then, there exists an orthonormal frame field $\{E_1,\cdot\cdot\cdot,E_n\}$ on $M^n$ satisfying:
\begin{align}\label{3.3}
\left \{
\begin{aligned}
& A_{E_1}E_1=\lambda_1 E_1,\quad A_{E_1}E_i=\mu E_i,\;\; \\
& A_{E_i}E_j=\mu\delta_{ij} E_1,\;\;\;\;\mu=\tfrac{\lambda_1}{2}>0,\quad 2\leq i,j\leq n.
\end{aligned}
\right.
\end{align}
if and only if the function $\theta$ defined by (\ref{3.2}) satisfies $\theta\equiv \lambda_1$, and the Calabi hypersurface $x$ has constant sectional curvature of $-\mu^2$, where $\mu=\frac{\lambda_1}{2}$.
\end{prop}

\textbf{Proof.}
For the ``only if'' part, if there exists an orthonormal frame field $\{E_1,\cdot\cdot\cdot,E_n\}$ on $M^n$ such that the Fubini-Pick tensor $A$ satisfies (\ref{3.3}), then it is obvious that the function $\theta =\lambda_1$.

Using the Codazzi equation (\ref{2.8}), we have
\begin{align}\label{3.4}
0=\;&\hat{\nabla}_{E_i}A(E_j,E_k)-A(\hat{\nabla}_{E_i}E_j,E_k)-
A(E_j,\hat{\nabla}_{E_i}E_k)-\hat{\nabla}_{E_j}A(E_i,E_k)\nnm\\
&+A(\hat{\nabla}_{E_j}E_i,E_k)+A(E_i,\hat{\nabla}_{E_j}E_k).
\end{align}
We write
$$\hat{\nabla}_{E_i}E_j= \hat{\Gamma}_{ij}^k E_k. $$
Choosing $i=1,\;2\leq j=k\leq n$  in (\ref{3.4}), we obtain
\begin{align}
0=\;& E_1(\mu)E_1+\mu\sum_{\alpha=1}^n\hat{\Gamma}_{11}^\alpha E_\alpha-2\mu\hat{\Gamma}_{1j}^1 E_j-E_j(\mu)E_j\notag\\ &-\mu\sum_{\alpha=1}^n\hat{\Gamma}_{jj}^\alpha E_\alpha+\mu\hat{\Gamma}_{j1}^j E_1+2\mu\hat{\Gamma}_{jj}^1E_1+\mu\sum_{l=2}^n\hat{\Gamma}_{jj}^l E_l.
\end{align}
Hence, we get
\begin{align*}
\begin{cases}
E_1(\mu)-\mu\hat{\Gamma}_{jj}^1+\mu\hat{\Gamma}_{j1}^j+2\mu\hat{\Gamma}_{jj}^1=0,\\
\mu\hat{\Gamma}_{11}^j-2\mu\hat{\Gamma}_{1j}^1-E_j(\mu)=0.\\
\end{cases}
\end{align*}
As $\hat{\Gamma}_{ij}^k+\hat{\Gamma}_{ik}^j=0$, we have
\begin{align}\label{3.6}
E_1(\mu)=0,\quad E_j(\mu)=3\mu\hat{\Gamma}_{11}^j,\;\ 2\leq j\leq n.
\end{align}
Similarly, choosing $i\geq 2,\;\;j=k=1$  in (\ref{3.4}), we obtain
\be\label{3.7} E_i(\mu)=0,\;\;  2\leq i\leq n.\ee
Then we have $\mu$ is a positive constant.

Taking in (\ref{2.7}) $i=k=1, j=l\geq 2$, we can get
$$R_{1j1j}=\sum_{m=1}^n (A_{1jm}^2-A_{11m}A_{jjm})=-\mu^2.$$
Taking in (\ref{2.7}) $i=k\geq 2, j=l\geq 2$ and $i\neq j$, we can get
$$R_{ijij}=\sum_{m=1}^n (A_{ijm}^2-A_{iim}A_{jjm})=-\mu^2.$$
Thus, $x$ is a Calabi hypersurface with constant sectional curvature $-\mu^2$.

Conversely, for the ``if'' part, we assume that $\theta\equiv \lambda_1>0$ and the sectional curvature $C=-\mu^2$, where $\mu=\frac{\lambda_1}{2}$ is a constant.

For any point $p\in M$, let $\{e_1,\cdots,e_n\}$ be the orthonormal basis of $T_{p}M^{n}$ in Lemma 3.1.  Using (\ref{3.1}), (\ref{2.7}) and $C=-\mu^2$,
we obtain that $\lambda_i=\frac{1}{2}\lambda_{1}=\mu,\,2\leq i\leq n$.
Hence, we have
$$A_{e_1}e_1=\lambda_{1}e_1,\;\; A_{e_1}e_i=\mu e_i,\;\; A_{e_i}e_j=\mu\delta_{ij}e_1, \;\; 2 \leq i,j \leq n. $$
Then the
Tchebychev vector field satisfies
$$T|_p=\tfrac{1}{n}\sum_{i=1}^{n} A(e_i,e_i)=\tfrac{1}{n}(\lambda_1+(n-1)\mu)e_1=\tfrac{n+1}{n}\mu e_1.$$
Put $E_1:=\frac {n}{(n+1)\mu}T$, then $E_1|_p=e_1$. As the function $\theta\equiv \lambda_1$, we see that $A_{E_1}E_1=\lambda_1 E_1.$ This implies that the distribution $\{E_1\}^\perp$ is $A_{E_1}$-invariant. Let $\{E_2,\cdots,E_n\}\in\{E_1\}^\perp$ be the orthonormal eigenvectors of $A_{E_1},$ then we have
$$A_{E_1}E_1= \lambda_1E_1,\; \; A_{E_1}E_i=\mu E_i, \;\;
A_{E_i}E_j=\mu\delta_{ij} E_1,\;\;\lambda_1=2\mu,\;\; 2\leq i,j\leq n.
\qquad\Box $$
\begin{rmk}
In the equiaffine differential geometry, this kind of frame field cannot occur because it has the apolarity condition. In the centroaffine geometry, this kind of frame field occurs and the corresponding centroaffine hypersurface is a canonical one, see Theorem 3.1 in \cite{ZCH}.
\end{rmk}

%%%%%%%%%%%%%%%%%%%%%%%%%%
%%%%%%%%%%%%%%%%%%%%%%%%%%
In the following, we assume that $M^n$ is an $n$-dimensional $(n\geq 3)$ Calabi hypersurface with constant sectional curvature $-\mu^2$, and the function $\theta \equiv \lambda_1=2\mu>0$ is a constant. Without loss of generality, we write $\mu=1$. By Proposition 3.2, there exists an orthonormal frame field $\left\{E_{1},\cdots,E_{n}\right\}$ on $M^n$, such that the Fubini-Pick tensor $A$ satisfies
\begin{align}\label{3.8}
\left \{
\begin{aligned}
&A_{E_1}E_1=2 E_1,\quad A_{E_1}E_i= E_i,\\
&A_{E_i} E_j=\delta_{ij} E_1, \quad 2 \leq i,j \leq n.
\end{aligned}
\right.
\end{align}

By the similar calculation in section 4 of \cite{CH}, we have the following lemmas 3.3-3.4.

\begin{lem}\label{3.3}
Let $\left\{E_{1},\cdots,E_{n}\right\}$ be the local orthonormal frame field on $M^n $ satisfying (\ref{3.8}), there hold:

\noindent{$\rm(i)$}\quad $\hat{\nabla}_{E_1}E_1=0;$

\noindent{$\rm(ii)$}\quad $\hat{\nabla}_{E_i}E_1=\eta E_i$, where $\eta:=-\hat{\Gamma}_{ii}^1,\  E_i(\eta)=0,\ 2\leq i\leq n;$

\noindent{$\rm(iii)$}\quad $E_1(\eta)+\eta^2-1=0.$
\end{lem}

\textbf{Proof.}
From (\ref{3.6}) and (\ref{3.7}), we obtain $\hat\nabla_{E_1}E_1=0$.
Choosing $2\leq i\leq n,2\leq j=k\leq n$ and $i\neq j$  in (\ref{3.4}), we have
$$
0=\sum_{\alpha=1}^n\hat{\Gamma}_{i1}^\alpha E_\alpha-2\hat{\Gamma}_{ij}^1 E_j+\hat{\Gamma}_{ji}^1 E_j+\hat{\Gamma}_{ji}^j E_1+\hat{\Gamma}_{jj}^1E_i+\hat{\Gamma}_{jj}^i E_1.
$$
Hence, we get
\begin{equation}\label{3.9}
\hat{\Gamma}_{ii}^1=\hat{\Gamma}_{jj}^1,\quad 2\leq i\neq j\leq n.
\end{equation}
Choosing distinct $i,j,k$ and $2\leq i,j,k\leq n$ in (\ref{3.4}), we have
$$
0=-\hat{\Gamma}_{ij}^1 E_k-\hat{\Gamma}_{ij}^k E_1-\hat{\Gamma}_{ik}^1 E_j-\hat{\Gamma}_{ik}^j E_1+\hat{\Gamma}_{ji}^1 E_k+\hat{\Gamma}_{ji}^k E_1+\hat{\Gamma}_{jk}^1E_i+\hat{\Gamma}_{jk}^i E_1.
$$
Hence, we get
\begin{equation}\label{3.10}
\hat{\Gamma}_{ik}^1=0,\quad i\neq k.
\end{equation}
From (\ref{3.9}) and (\ref{3.10}), we have
$$\hat{\nabla}_{E_i}E_1=\hat{\Gamma}_{i1}^\alpha E_\alpha=-\hat{\Gamma}_{i\alpha}^1 E_\alpha=-\hat{\Gamma}_{ii}^1 E_i:=\eta E_i,\quad 2\leq i\leq n.$$

By the definition of the Riemannian curvature tensor $\mathcal{R}$ of the Calabi metric and the Gauss equation, we have
\begin{align}\label{3.11}
\mathcal{R}(E_i,E_j)E_k=&A_{E_j}A_{E_i}E_k-A_{E_i}A_{E_j}E_k\notag\\
=&\hat{\nabla}_{E_i}\hat{\nabla}_{E_j}E_k-\hat{\nabla}_{E_j}\hat{\nabla}_{E_i}E_k-\hat{\nabla}_{[E_i,E_j]}E_k.
\end{align}
Choosing $k=1,2\leq i\neq j\leq n$ in (\ref{3.11}), we have
$$0=A_{E_j}E_i-A_{E_i}E_j=\hat{\nabla}_{E_i}(\eta E_j)-\hat{\nabla}_{E_j}(\eta E_i)-\eta [E_i,E_j].$$
Then we get $E_i(\eta)E_j-E_j(\eta)E_i=0$ for $i\neq j$. It implies that $E_i(\eta)=0$ for $2\leq i\leq n$.
Choosing $i=k=1,2\leq j\leq n$ in (\ref{3.11}), we have
\begin{align}\label{3.12}
2E_j-E_j=&\hat{\nabla}_{E_1}(\eta E_j)-\hat{\nabla}_{E_j}(\hat{\nabla}_{E_1}E_1)-\eta [E_1,E_j]\notag\\
=&E_1(\eta)E_j+\eta\hat{\nabla}_{E_1}E_j-\eta(\hat{\nabla}_{E_1}E_j-\hat{\nabla}_{E_j}E_1)\notag\\
=&E_1(\eta)E_j+\eta^2E_j.
\end{align}
Then we get $E_1(\eta)+\eta^2-1=0.$   \hfill $\square$

\begin{rmk} In the two dimension case, Lemma 3.3 remains valid with the exception that $E_i(\eta)=0$ for $i=2$ in (ii).
However, the lemma retains its validity in full when the Calabi surface satisfies the affine maximal type equation.
\end{rmk}
To state the following Lemma 3.4, we write $\mathcal{D}:=span\{E_2,\cdot\cdot\cdot,E_n\}$.

\begin{lem}
Around any $p\in M^n$, there hold:

\noindent{$\rm(i)$}\; $\eta$ can be expressed as a function $\eta(t)$ defined in an open interval $I$ containing 0, and it satisfies the real ordinary differential equation
\be\label{3.13}
\eta'+\eta^2-1=0.\ee
Moreover, writing $\rho:=\rho(t)=e^{\int \eta dt}$, the function $\rho^2(\eta^2-1)$ is a constant and we denote it by $\bar{c}$;

\noindent{$\rm(ii)$}\;
$(M^n,G)$ is locally isometric to the warped product $I\times_{\rho}M_1=(I\times M_1,dt^2+\rho^2 G_1)$, where, $M_1$ is an integral manifold of the distribution $\mathcal{D}$, and $(M_1,G_1)$ is a space form with constant sectional curvature $\bar{c}$.
\end{lem}

\textbf{Proof.}
Let $\left\{E_{1},...,E_{n}\right\}$ be the local orthonormal frame field on a neighbourhood of $p\in M^n$ satisfying (\ref{3.8}). The fact $(i)$ of Lemma 3.3 implies that the distribution $\mathcal{D}^{\perp}$ is auto-parallel. Then, by $(ii)$ of Lemma 3.3, for any $E_i,E_j\in\mathcal{D},2\leq i,j\leq n$, we have that
\be G(\hat{\nabla}_{E_i}E_j,E_1)=-G(E_j,\hat{\nabla}_{E_i}E_1)=-\eta G(E_i,E_j)=G(-\eta E_1,E_1)G(E_i,E_j),\ee
and
\be G(\hat{\nabla}_{E_i}(-\eta E_1),E_1)=0.\ee
Hence, $\mathcal{D}$ is spherical. Therefore, by applying Theorem 2.1, $(M^n,G)$ is locally isometric to a warped product $I\times_\rho M_1$ where $I$ is an open interval containing $0$, $M_1$ is the integral manifold of the distribution $\mathcal{D}$ and $\rho$ is the warping function which is determined by
\begin{equation}\label{3.16}
-\eta E_1=-E_1(\ln\rho)E_1.
\end{equation}

Next, the warped product structure of $(M^n,G)$ implies the existence of local coordinates $(t,u)\in \mathbb{R}\times \mathbb{R}^{n-1}$ for $M^n$ with coordinates domain around the point $p$ (corresponding to the origin in $\mathbb{R}^n$), such that $\mathcal D^\perp$ is given by $du=0$ and $\mathcal D$ is given by $dt=0$. Moreover, we may assume that $E_1=\frac\partial{\partial t}$, and then $G=dt^2+\rho^2G_1$, where $G_1$ is the metric tensor of $M_1$.
From $(ii)$ in Lemma 3.3, we know that $\eta$ depends only on $t$, and we write $\eta=\eta(t)$. Then $(iii)$ of
Lemma 3.3 imply that (\ref{3.13}) holds.

As a warped product, it is well known that the curvature tensor $\mathcal{R}$ of $(M^n,G)$ is related to the curvature tensor $\mathcal{R}^{M_1}$ of $(M_1,G_1)$ by
\be \label{3.17}
\mathcal{R}(X,Y)Z=\mathcal{R}^{M_1}(X,Y)Z-(\tfrac{\rho^{'}}{\rho})^2(G(Y,Z)X-G(X,Z)Y).
\ee
By (\ref{3.8}), (\ref{3.16}) and (\ref{3.17}), for any $2\leq i,j,k,l\leq n$, we have
\begin{align} \label{3.18}
G_1(\mathcal{R}^{M_1}(E_i,E_j)E_k,E_l)=&\;G_1(A_{E_j}A_{E_i}E_k-A_{E_i}A_{E_j}E_k+\eta^{2}
(\delta_{jk}E_i-\delta_{ik}E_j),E_l)\notag\\
=&\tfrac{\eta^2-1}{\rho^{2}}(\delta_{jk}\delta_{il}-\delta_{ik}\delta_{jl}).
\end{align}
Then the sectional curvature of $(M_1,G_1)$ is
\begin{equation}\label{3.19}
K^{M_1}(E_i,E_j)
=\rho^{2}(\eta^2-1).
\end{equation}
Moreover, by using (\ref{3.13}) and (\ref{3.16}), a direct calculation shows that $\rho^{2}(\eta^2-1)$ is actually a constant. This fact and (\ref{3.19}) imply that $(M_1,G_1)$ is a space form with constant sectional curvature $\bar{c}:=\rho^{2}(\eta^2-1).$ \hfill $\square$

\section{Classification of Calabi hypersurfaces with such frame}

Having the above Lemma 3.4 and using the similar method in \cite{LLX}, we will prove the following classification result.

\begin{thm}
Let $n\geq3$ and $M^n$ be a Calabi hypersurface with constant sectional curvature $-1$ in $\mathbb{R}^{n+1}$. If the function $\theta\equiv 2$, then, corresponding to each of the three possibilities depending on the sign of $\bar{c}$, the hypersurface $M^n$ is Calabi affine equivalent to an open part of one of the following hypersurfaces:
\begin{enumerate}
  \item Case $\bar{c}>0$,\begin{equation}\label{4.1}
x_{n+1}=\tfrac{1}{4}(\sum_{k=1}^nx_k^2)^{\tfrac{1}{2}}-\tfrac{1}{4}\ln{[(\sum_{k=1}^nx_k^2)^{\tfrac{1}{2}}+1]},
\quad\sum_{k=1}^nx_k^2>0.
\end{equation}
  \item  Case $\bar{c}=0$,\begin{equation}\label{4.2}
x_{n+1}=-\tfrac{1}{4}\ln{(x_1-\sum^n_{k=2}\tfrac{x^2_k}{2})},\;\;\;\;\;  x_1>\sum^n_{k=2}\tfrac{x^2_k}{2};
\end{equation}
or
\begin{equation}\label{4.3}
x_{n+1}=-\tfrac{1}{4}\ln{x_1}+\sum^n_{k=2}\tfrac{x^2_k}{2x_1},\;\;\;\;\;  x_1>0.
\end{equation}
  \item Case $\bar{c}<0$,\begin{equation}\label{4.4}
x_{n+1}=-\tfrac{1}{4}(x_1^2-\sum_{k=2}^nx_k^2)^{\frac{1}{2}}-\tfrac{1}{4}
\ln{[(x_1^2-\sum_{k=2}^nx_k^2)^{\frac{1}{2}}-1]},
\;\; x_1^2-\sum_{k=2}^nx_k^2>1.
\end{equation}
\end{enumerate}
\end{thm}

\textbf{Proof.} The discussions for all these three subcases are very
similar, we will give the details for Case (1) and the main formulas for Case (2) and Case (3).

\textbf{Case (1).} In this case, by Lemma 3.4, $(M^n,G)$ is locally isometric to a warped product $I\times_{\rho}\mathbb{S}^{n-1}(1/\sqrt{\bar{c}})$. Up to a Calabi affine equivalent, we may assume $\bar{c}=1$. From (\ref{3.13}), (\ref{3.16}) and $\bar{c}=\rho^2(\eta^2-1)$, we obtain that $\rho=\pm\sinh{(t+c_1)}$, where $c_1$ is a constant. Without loss of generality, we may assume that $\rho=\sinh{t}$, where $t>0$. Hence, we may choose local coordinates $(t,u_2,\cdot\cdot\cdot,u_n)$ on $M^n$ such that the metric $G$ has the following expression:
\be\label{4.5}
G=dt^2+\sinh^2{t}\bigg(du_2^2+\sin^2{u_2}du_3^2+
\cdot\cdot\cdot+\prod^{n-1}_{i=2}
\sin^2{u_i} du_n^2\bigg).\ee

From (\ref{4.5}), we can calculate the Levi-Civita connections as follows:
\begin{equation}\label{4.6}
\left \{
\begin{aligned}
 \hat{\nabla}_{\partial t}\partial t=0,&\; \hat{\nabla}_{\partial u_2}\partial u_2=-\sinh{t}\cosh{t}\partial t,\; \hat{\nabla}_{\partial t}\partial u_k=\coth{t}\partial u_k,\; 2\leq k\leq n,\\
 \hat{\nabla}_{\partial u_k}\partial u_j
=&\;\cot{u_k}\partial u_j,\ \;\; 2\leq k<j\leq n ,\\
 \hat{\nabla}_{\partial u_3}\partial u_3
=&-\sinh{t}\cosh{t}\sin^2{u_2}\partial t-\tfrac{1}{2}\sin(2u_2)\partial u_2,\\
 \hat{\nabla}_{\partial u_k}\partial u_k
=&-\sinh{t}\cosh{t}\prod^{k-1}_{j=2}\sin^2 u_j\partial t\\
 &-\sum_{j=2}^{k-1}\tfrac{1}{2}\sin(2u_j)\prod^{k-1}_{i=j+1}
\sin^2 u_i\partial u_j, 4\leq k\leq n,
\end{aligned}
\right.
\end{equation}
where, and later, we use the notations $\partial t=\frac{\partial}{\partial t},\ \partial u_k=\frac{\partial}{\partial u_k}$.

By using (\ref{3.8}) and (\ref{4.5}), we get that
$$A_{\partial t}X=X,\ \;\; A_XZ=G(X,Z)\partial t,\quad  X,Z\in\mathcal{D}.$$
It follows that
\begin{equation}\label{4.7}
\left \{
\begin{aligned}
&A_{\partial t}\partial t=2\partial t,\;\; A_{\partial t}\partial u_k=\partial u_k,\ \;\; 2\leq k\leq n, \\
&A_{\partial u_2}\partial u_2=\sinh^2{t}\partial t,\;\; A_{\partial u_k}\partial u_j=0,\ \;\; 2\leq j,k\leq n \;\; \text{and}\;\;  j\neq k,\\
&A_{\partial u_k}\partial u_k=\sinh^2{t}\prod^{k-1}_{j=2}
\sin^2 u_j\partial t,\ \;\; 3\leq k\leq n.
\end{aligned}
\right.
\end{equation}
Then, using the definition of $A$ and a straightforward computation, we get the Gauss equation for $x$,
\begin{align}
\label{4.8} x_{tt}=&\; 2x_{t}+Y, \\
\label{4.9}x_{ tu_k}=&\;(1+\coth{t})x_{ u_k},\qquad\qquad 2\leq k\leq n,\\
\label{4.10} x_{u_2 u_2}=&\;\sinh^2{t}(1-\coth{t}) x_{ t}+\sinh^2{t}Y, \\
\label{4.11}x_{u_ku_j}=&\;\cot u_kx_{u_j},\qquad\qquad 2\leq k<j\leq n,\\
\label{4.12}x_{u_3u_3}=&\;\sinh^2{t}(1-\coth{t})\sin^2 u_2x_t
-\tfrac{1}{2}\sin (2u_2)x_{u_2}+\sinh^2{t}\sin^2 u_2Y,
\end{align}
\begin{align}
\label{4.13}x_{u_ku_k}=&\;\sinh^2{t}(1-\coth{t})\prod^{k-1}_{j=2}\sin^2u_jx_t-\sum^{k-1}_{j=2}
\bigg(\tfrac{1}{2}\sin(2u_j)\prod^{k-1}_{i=j+1}\sin^2u_i\bigg)x_{u_j}\notag\\
&+\sinh^2{t}\prod^{k-1}_{j=2}\sin^2u_jY,\qquad\qquad 4\leq k\leq n,
\end{align}
where $Y=(0,\cdots,0,1)$. To solve the above equations, we first have to solve their corresponding homogeneous equations.
\begin{align}
\label{4.14} x_{tt}=&\; 2x_{t}, \\
\label{4.15} x_{ tu_k}=&\;(1+\coth{t})x_{ u_k},\qquad  2\leq k\leq n, \\
\label{4.16} x_{u_2 u_2}=&\;\sinh^2{t}(1-\coth{t}) x_{ t}, \\
\label{4.17} x_{u_ku_j}=&\;\cot u_kx_{u_j},\qquad\qquad 2\leq k<j\leq n,\\
\label{4.18} x_{u_3u_3}=&\;\sinh^2{t}(1-\coth{t})\sin^2 u_2x_t-\tfrac{1}{2}\sin (2u_2)x_{u_2}, \\
\label{4.19} x_{u_ku_k}=&\;\sinh^2{t}(1-\coth{t})\prod^{k-1}_{j=2}\sin^2u_jx_t\notag\\
 &\; -\tfrac{1}{2}\sum^{k-1}_{j=2}\sin(2u_j)\prod^{k-1}_{i=j+1}
\sin^2u_ix_{u_j}, \qquad 4\leq k\leq n.
\end{align}
We can solve (\ref{4.14}) to obtain that
\begin{align}\label{4.20}
x_0 = Q_1(u_2,\cdot\cdot\cdot,u_n)\int e^{2t}
+Q_2(u_2,\cdot\cdot\cdot,u_n),
\end{align}
where $Q_1(u_2,\cdot\cdot\cdot,u_n)$ and $Q_2(u_2,\cdot\cdot\cdot,u_n)$ are $\mathbb{R}^{n+1}$-valued functions. Inserting (\ref{4.20}) into (\ref{4.15}), we obtain
\be\label{4.21}
\tfrac{\partial Q_1}{\partial u_k}\left(\tfrac{e^{2t}}{1+\coth{t}}-\int e^{2t}\right)=\tfrac{\partial Q_2}{\partial u_k}.\ee

On the other hand, we can easily check that
$$\tfrac{e^{2t}}{1+\coth{t}}-\int e^{2t}=const.$$
Without loss of generality, we may choose the antiderivative of $e^{2t}$ such that
\be\label{4.22}
\tfrac{e^{2t}}{1+\coth{t}}-\int e^{2t}=0.\ee
Then we obtain
$$\int e^{2t}=\tfrac{1}{2}e^{2t}-\tfrac{1}{2}.$$

From (\ref{4.21}) and (\ref{4.22}), we obtain $\frac{\partial Q_2}{\partial u_k}=0$ for $2\leq k\leq n,$ which implies that $Q_2(u_2,\cdot\cdot\cdot,u_n)=A_1$ is a constant vector. Hence, we have
\be\label{4.23}
x_0 = Q_1(u_2,\cdot\cdot\cdot,u_n)(\tfrac{1}{2}e^{2t}-\tfrac{1}{2})+A_1.\ee
Putting (\ref{4.23}) into (\ref{4.16}), and noting that
\be \label{4.24} \sinh^2t (1-\coth t)=\tfrac{-1}{1+\coth t}=\tfrac{1}{2}(e^{-2t}-1),\ee
we obtain that $\frac{\partial^2 Q_1}{\partial u_2 \partial u_2}=-Q_1$. Then we have
\be\label{4.25}
Q_1=Q_3(u_3,\cdot\cdot\cdot,u_n)\cos u_2+Q_4(u_3,\cdot\cdot\cdot,u_n)\sin u_2.\ee
From (\ref{4.25}) and (\ref{4.17}) with $k=2$, we conclude that $Q_3(u_3,\cdot\cdot\cdot,u_n)=A_2$ is a constant vector. Thus, together with (\ref{4.25}), (\ref{4.18}) and (\ref{4.24}), we get
\be\label{4.26}
Q_4(u_3,\cdot\cdot\cdot,u_n)=Q_5(u_4,\cdot\cdot\cdot,u_n)\cos u_3
+Q_6(u_4,\cdot\cdot\cdot,u_n)\sin u_3.\ee

Similarly, by solving (\ref{4.26}) and (\ref{4.17}) with $k=3$, we conclude that $Q_5=A_3$ is a constant vector. Carrying such procedure by induction, we finally obtain that:
\begin{align}\label{4.27}
x_0=&A_1+ (\tfrac{1}{2}e^{2t}-\tfrac{1}{2})\bigg(A_2\cos u_2+A_3\sin u_2\cos u_3+A_4\sin u_2\sin u_3\cos u_4+\cdot\cdot\cdot\notag\\
&+A_n\prod_{i=2}^{n-1}\sin u_i\cos u_n+A_{n+1}\prod_{i=2}^{n}\sin u_i\bigg),
\end{align}
where $A_k(1\leq k\leq n+1)$ are all constant vectors. On the other hand, we know that
\be\label{4.28}
\bar{x}=\bigg(0,\cdot\cdot\cdot,0,\tfrac{1}{4}e^{2t}-\tfrac{t}{2}\bigg)\ee
is a special solution of equations (\ref{4.8})-(\ref{4.13}).
Consequently, the general solution of equations (\ref{4.8})-(\ref{4.13}) can be expressed as
\begin{align} \label{4.29}
x=&A_1+(\tfrac{1}{2}e^{2t}-\tfrac{1}{2})\bigg(A_2\cos u_2+A_3\sin u_2\cos u_3+A_4\sin u_2\sin u_3\cos u_4+\cdot\cdot\cdot\notag\\
&+A_n\prod_{i=2}^{n-1}\sin u_i\cos u_n+A_{n+1}\prod_{i=2}^{n}\sin u_i\bigg)+\bigg(0,\cdot\cdot\cdot,0,\tfrac{1}{4}e^{2t}-\tfrac{t}{2}\bigg).
\end{align}

Since $M^n$ is non-degenerate, $x-A_1$ lies linearly full in $\mathbb{R}^{n+1}$. Then $A_2,\ \cdot\cdot\cdot,\ A_{n+1}$ and $(0, \cdot\cdot\cdot, 0,1)$ are linearly independent vectors. Thus there exists an affine transformation $\phi\in SA(n+1)$ such that
\begin{align*}
&A_1=(0,0,\cdot\cdot\cdot,0),\ A_2=(1,0,\cdot\cdot\cdot,0),\  A_3=(0,1,0,\cdot\cdot\cdot,0),\  \cdot\cdot\cdot,\\
& A_n=(0,\cdot\cdot\cdot,0,1,0,0),\ A_{n+1}=(0,\cdot\cdot\cdot,0,1,0).
\end{align*}
Then, the position vector $x$ can be written as
\begin{equation}\label{4.30}
\left \{
\begin{aligned}
x_1=\;&(\tfrac{1}{2}e^{2t}-\tfrac{1}{2})\cos u_2, \\
x_k=&(\tfrac{1}{2}e^{2t}-\tfrac{1}{2})\prod_{i=2}^{k}\sin u_i\cos u_{k+1},\  2\leq k\leq n-1,\\
x_n=&(\tfrac{1}{2}e^{2t}-\tfrac{1}{2})\prod_{i=2}^{n}\sin u_i,\qquad
x_{n+1}=\;  \tfrac{1}{4}e^{2t}-\tfrac{t}{2}.
\end{aligned}
\right.
\end{equation}
Therefore, up to an affine transformation, $M^n$ locally lies on the graph hypersurface of the function (\ref{4.1}).

\textbf{Case (2).} When $\bar{c}=0$, $(M^n,G)$ is locally isometric to a warped product $I\times_{\rho}\mathbb{R}^{n-1}$. Hence, we may choose local coordinates $(t,u_2,\cdot\cdot\cdot,u_n)$ on $M^n$ such that the metric $G$ has the following expression:
\be\label{4.31}
G=dt^2+\rho^2(du_2^2+du_3^2+\cdot\cdot\cdot+du_n^2).\ee
It follows from (\ref{4.31})  that
\begin{equation}\label{4.32}
\left \{
\begin{aligned}
& \hat{\nabla}_{\partial t}\partial t=0,\   \hat{\nabla}_{\partial t}\partial u_k=\frac{\rho'}{\rho}\partial u_k,\quad 2\leq k\leq n, \\
&  \hat{\nabla}_{\partial u_k}\partial u_j=0,\ 2\leq j,k\leq n, \ \,\qquad  j\neq k,\\
&  \hat{\nabla}_{\partial u_k}\partial u_k=-\rho'\rho\partial t, \qquad \qquad2\leq k\leq n.\\
\end{aligned}
\right.
\end{equation}
By using (\ref{3.8}) and (\ref{4.31}), we get that
\begin{equation}\label{4.33}
\left \{
\begin{aligned}
&A_{\partial t}\partial t=2\partial t,\ A_{\partial t}\partial u_k=\partial u_k,\qquad 2\leq k\leq n, \\
&A_{\partial u_k}\partial u_j=0,\ 2\leq j,k\leq n, \ \,\qquad  j\neq k,\\
&A_{\partial u_k}\partial u_k=\rho^2\partial t,\qquad \qquad 2\leq k\leq n.\\
\end{aligned}
\right.
\end{equation}

According to the expression for $\bar{c}$ and $\rho \neq 0$, we get $\eta=-1$ or $\eta=1$. The following two cases are calculated separately.

 \textbf{Case (2.1).} $\eta=-1$.

We can solve (\ref{3.16}) to obtain that $\rho=c_2e^{-t}$, where $c_2$ is a positive constant.
Up to a Calabi affine equivalent, we may assume $c_2=1$. Then, from (\ref{4.32}),(\ref{4.33}) and using the definition of $A$, we get the Gauss equation for $x$
\begin{align}
\label{4.34} x_{tt}&=2 x_t+Y, \\
\label{4.35} x_{tu_k}&=0,\qquad\qquad \qquad\qquad 2\leq k\leq n, \\
\label{4.36} x_{u_ku_k}&=2e^{-2t}x_t+e^{-2t}Y,\qquad 2\leq k\leq n, \\
\label{4.37} x_{u_ku_j}&=0,\qquad\qquad \qquad 2\leq j,k\leq n,  \;  \ j\neq k.
\end{align}
As in \textbf{Case(1)}, the general solution of equations (\ref{4.34})-(\ref{4.37}) can be expressed as
\be \label{4.38}
x= A_2\bigg(2\sum^n_{k=2}u^2_k+e^{2t}\bigg)+\sum^n_{k=2}A_{k+1}u_k+A_{1}
+\bigg(0,\cdot\cdot\cdot,0,-\tfrac{1}{2}t\bigg),
\ee
where $A_1,\ \cdot\cdot\cdot,\ A_{n+1}$ are constant vectors in $\mathbb {R}^{n+1}$.
Then, the position vector $x$ can be written as
\be \label{4.39}
 x_1=2\sum^n_{k=2}u^2_k+e^{2t}; \quad
 x_k=u_k,\quad 2\leq k\leq n;\quad
 x_{n+1}=-\tfrac{1}{2}t.
\ee
Therefore, up to an affine transformation, $M^n$ locally lies on the graph hypersurface of the function (\ref{4.2}).

\noindent\textbf{Case (2.2).} $\eta=1.$

Similarly, using the definition of $A$ and a straightforward computation, we get the Gauss equation for $x$,
\begin{align}
\label{4.40}x_{tt}&=2 x_t+Y, \\
\label{4.41}x_{tu_k}&=2 x_{u_k},\;\qquad 2\leq k\leq n, \\
\label{4.42}x_{u_ku_k}&=e^{2t}Y,\;\qquad 2\leq k\leq n, \\
\label{4.43}x_{u_ku_j}&=0,\ \qquad  2\leq j,\;  k\leq n,\;\;   j\neq k.
\end{align}
By the same method as in \textbf{Case(1)}, we obtain that the general solution of equations (\ref{4.40})-(\ref{4.43}) is
\begin{align}\label{4.44}
x=A_{1}+\bigg(\sum^n_{k=2}A_{k+1}u_k+A_{2}\bigg)e^{2t}+\left(0,\cdot\cdot\cdot,0,\tfrac{1}{2}e^{2t}
\sum^n_{k=2}u_k^2-\tfrac{t}{2}\right),
\end{align}
where $A_1,\ \cdot\cdot\cdot,\ A_{n+1}$ are constant vectors in $\mathbb {R}^{n+1}$.
Then, the position vector $x$ can be written as
\be \label{4.45}
x_1=e^{2t};\;\;
x_k=e^{2t}u_k,\ \;\; 2\leq k\leq n;\;\;
x_{n+1}=\tfrac{1}{2}e^{2t}\sum^n_{k=2}u_k^2-\tfrac{t}{2}.
\ee
Therefore, $M^n$ locally lies on the graph hypersurface of the function (\ref{4.3}).

\noindent\textbf{Case (3).} In this case, $(M^n,G)$ is locally isometric to a warped product $I\times_{\rho}\mathbb{H}^{n-1}(1/\sqrt{-\bar{c}})$.
Up to a Calabi affine equivalent, we may assume $\bar{c}=-1$. From (\ref{3.13}), (\ref{3.16}) and  $\bar{c}=\rho^2(\eta^2-1)$, we obtain that $\rho=\pm\cosh{(t+c_3)}$, where $c_3$ is a constant.
Without loss of generality, we assume $\rho=\cosh{t} $. Hence, we can choose local coordinates $(t,u_2,\cdot\cdot\cdot,u_n)$ on $M^n$ such that the metric $G$ has the following expression:
\begin{equation}\label{4.46}
\begin{aligned}
G=&dt^2+\cosh^2{t}\big[du_2^2+\sinh^2 u_2\big(du_3^2+\sin^2 u_3du^2_4 \\ &+\cdot\cdot\cdot+\sin^2 u_3\cdot\cdot\cdot
\sin^2 u_{n-1}du_n^2\big)\big].
\end{aligned}
\end{equation}

From (\ref{4.46}), we can calculate the Levi-Civita connections as follows:
\begin{equation}\label{4.47}
\left \{
\begin{aligned}
\hat{\nabla}_{\partial t}\partial t=0,& \hat{\nabla}_{\partial t}\partial u_k=\tanh{t}\partial u_k,\; 2\leq k\leq n, \; \hat{\nabla}_{\partial u_2}\partial u_2=-\sinh{t}\cosh{t}\partial t,\\
\hat{\nabla}_{\partial u_2}\partial u_k=&\coth u_2\partial u_k,\;\; 3\leq k\leq n,\\
\hat{\nabla}_{\partial u_3}\partial u_3=&-\sinh{t}\cosh{t}\sinh^2 u_2\partial t-\tfrac{1}{2}\sinh(2u_2)\partial u_2,\\
\hat{\nabla}_{\partial u_k}\partial u_j=&\cot u_k\partial u_j,\;\; 3\leq k<j\leq n,\\
\hat{\nabla}_{\partial u_k}\partial u_k=&-\prod^{k-1}_{j=3}\sin^2 u_j(\sinh{t}\cosh{t}\sinh^2 u_2\partial t+\tfrac{1}{2}\sinh(2u_2)\partial u_2) \\
&-\sum_{j=3}^{k-1}\tfrac{1}{2} \sin(2u_j)\prod^{k-1}_{i=j+1}\sin^2 u_i\partial u_j ,\;\;\; 4\leq k\leq n.
\end{aligned}
\right.
\end{equation}
By using (\ref{3.8}) and (\ref{4.46}), we get
\begin{equation}\label{4.48}
\left \{
\begin{aligned}
&A_{\partial t}\partial t=2\partial t,\;\;\;\ A_{\partial t}\partial u_k=\partial u_k,\ 2\leq k\leq n, \\
&A_{\partial u_2}\partial u_2=\cosh^2{t}\partial t,
\;\;\;A_{\partial u_k}\partial u_j=0,\;\;\ 2\leq j,k\leq n, \ \;\;\  j\neq k,\\
&A_{\partial u_3}\partial u_3=\cosh^2{t}\sinh^2 u_2\partial t,\\
&A_{\partial u_k}\partial u_k=\cosh^2{t}\sinh^2 u_2\prod^{k-1}_{j=3}\sin^2 u_j\partial t,\;\;\; 4\leq k\leq n.
\end{aligned}
\right.
\end{equation}
Then, using the definition of $A$ and a straightforward computation, we get the Gauss equation for $x$,
\begin{align}
\label{4.49}x_{tt}=&\;2x_t+Y, \\
\label{4.50}x_{tu_k}=&\;(1+\tanh{t})x_{u_k},\ 2\leq k\leq n, \\
\label{4.51}x_{u_2u_2}=&\;\cosh^2{t}(1-\tanh{t})x_t+\cosh^2{t}Y, \\
\label{4.52}x_{u_2u_k}=&\;\coth u_2x_{u_k},\ 3\leq k\leq n,\\
\label{4.53}x_{u_3u_3}=&\;\cosh^2{t}(1-\tanh{t})\sinh^2 u_2x_t-\tfrac{1}{2}\sinh (2u_2)x_{u_2}+\cosh^2{t}\sinh^2 u_2Y,\\
\label{4.54}x_{u_ku_j}=&\;\cot u_k x_{u_j},\ 3\leq k<j\leq n,\\
\label{4.55}x_{u_ku_k}=&\;\cosh^2{t}(1-\tanh{t})\sinh^2 u_2\prod^{k-1}_{j=3}\sin^2 u_jx_t-\tfrac{1}{2}\sinh (2u_2)\prod^{k-1}_{j=3}\sin^2 u_jx_{u_2} \notag\\
-\tfrac{1}{2}\sum^{k-1}_{j=3}&\sin(2u_j)\prod^{k-1}_{i=j+1}\sin^2u_ix_{u_j}
 +\cosh^2{t}\sinh^2 u_2\prod^{k-1}_{j=3}\sin^2 u_jY,\ 4\leq k\leq n.
\end{align}

By the same method as in \textbf{Case(1)}, we can obtain that the general solution of equations (\ref{4.49})-(\ref{4.55}) is
\begin{align}\label{4.56}
x =&A_1+(\tfrac{1}{2}e^{2t}+\tfrac{1}{2})\bigg[A_2\cosh u_2+\sinh u_2\big(A_3\cos u_3+A_4\sin u_3\cos u_4+\cdot\cdot\cdot\notag\\
&+A_n\prod_{i=3}^{n-1}\sin u_i\cos u_n+A_{n+1}\prod_{i=3}^{n}\sin u_i\big)\bigg]
+\bigg(0,\cdot\cdot\cdot,0,-\tfrac{1}{4}e^{2t}-\tfrac{t}{2}\bigg).
\end{align}
Then, there exists an affine transformation $\phi\in SA(n+1)$ such that the position vector $x$ can be written as
\begin{equation}\label{4.57}
\left \{
\begin{aligned}
x_1=&\;(\tfrac{1}{2}e^{2t}+\tfrac{1}{2})\cosh u_2, \quad
x_2= \;(\tfrac{1}{2}e^{2t}+\tfrac{1}{2})\sinh u_2\cos u_3,\\
x_k=&\;(\tfrac{1}{2}e^{2t}+\tfrac{1}{2})\sinh u_2\prod_{i=3}^{k}\sin u_i\cos u_{k+1},\  3\leq k\leq n-1,\\
x_n=&\;(\tfrac{1}{2}e^{2t}+\tfrac{1}{2})\sinh u_2\prod_{i=3}^{n}\sin u_i, \quad
x_{n+1}= -\tfrac{1}{4}e^{2t}-\tfrac{t}{2}.\\
\end{aligned}
\right.
\end{equation}
Therefore, up to an affine transformation, $M^n$ locally lies on the graph hypersurface of the function (\ref{4.4}).\hfill $\Box$

\section{Proof of Theorem 1.3}

Define a domain
\begin{align*}
\Omega:=&\{(x_1,\dots,x_n)\;|\; (x_1,\dots,x_n)\in\mathbb{R}^n,\;x_1>\sum^n_{k=2}\tfrac{x^2_k}{2}\}.
\end{align*}
Consider the graph hypersurface $M^n:=\{ (x,f(x))\;|\;(x_1,\dots,x_n)\in \Omega \}$,
where
\be \label{5.1}
f(x):= \;-\tfrac{1}{4}\ln{(x_1-\sum^n_{k=2}\tfrac{x^2_k}{2})}.
\ee
Then it is easy to see that $M^n$ is a Euclidean complete Calabi hypersurface.

Denote $g(x):=x_1-\sum^n_{k=2}\frac{x^2_k}{2}>0$ on $\Omega$. Note that $f(x)=-\frac{1}{4}\ln g$, then we have
\begin{align*}
f_1=&\tfrac{\partial f}{\partial x_1}=-\tfrac{1}{4g},\ \;\;\quad f_k=\tfrac{x_k}{4g},\;\; k\geq 2,\;\;\quad f_{11}=\tfrac{1}{4g^2},\\
f_{1k}=& -\tfrac{x_k}{4g^2},\;\; f_{kk}=\tfrac{1}{4g^2}(g+x_k^2),\;\; f_{kj}=\tfrac{x_kx_j}{4g^2},\;\;2\leq k\neq j\leq n.
\end{align*}
The determinant and inverse matrix of the Hessian matrix $(f_{ij})$ are
\begin{equation}\label{5.2}
D:=\det{(f_{ij})}=\tfrac{1}{4^ng^{n+1}},
\end{equation}
\begin{equation}\label{5.3}
(f^{ij})=4g
\left(
\begin{array}{ccccc}
2x_1-g&x_2&x_3&\cdots&x_n\\
x_2&1&0&\cdots&0\\
x_3&0&1&\cdots&0\\
\vdots &\vdots &\vdots &\ddots &\vdots \\
x_n&0&0&\cdots&1
\end{array}
\right).
\end{equation}
From (\ref{5.2}) and (\ref{5.3}), we have
\begin{equation}\label{5.4}
\sum f^{ij}(D^a)_{ij}=4(n+1)[(n+1)a^2+na]D^{a},\quad a\neq 0.
\end{equation}
Hence the function $f(x)$ satisfies the affine maximal type equation (\ref{1.2}) if and only if $a=-\frac{n}{n+1}$.
If $a=-1$, then $f(x)$ satisfies the PDE (\ref{1.1}) with $L^\sharp=-4(n+1)$.

From (\ref{2.5}), the norm square of the Tchebychev vector field of $M^n$
\begin{equation}\label{5.5}
|T|^2=\tfrac{1}{4n^2}\sum f^{ij}(\ln{D})_{i}(\ln{D})_{j}=\tfrac{(n+1)^2}{n^2}.
\end{equation}
Then, by using (\ref{2.10}) and Theorem 4.1, we get the Pick invariant $J$
\be \label{5.6}
J=C+\tfrac{n}{n-1}|T|^2=\tfrac{3n+1}{n(n-1)}.
\ee
Thus, from Theorem 5.9.1 in \cite{LXSJ} or Theorem 2.10 in \cite{CLS}, we easily know that $M^n$ is complete with respect to the Calabi metric, see also Proposition 3.8 in \cite{S}.

Denote
\begin{align*}
\tilde{\Omega}:=\{(y_1,\dots,y_n)\;|\; (y_1,\dots,y_n)\in\mathbb{R}^n,\; y_1>0\}.
\end{align*}
Consider the graph hypersurface
$\tilde M^n:=\{(y,\tilde{f}(y))\;|\;(y_1,\dots,y_n)\in \tilde{\Omega}\}$,
where
\be \label{5.7}
\tilde{f}(y):= \;-\tfrac{1}{4}\ln{y_1}+\sum^n_{k=2}\tfrac{y^2_k}{2y_1}.\ee
Then it is easy to see that $\tilde M^n$ is also a Euclidean complete Calabi hypersurface.
In the sense of Calabi-affine equivalence, the function $f(x)$ given by (\ref{5.1}) is exact the Legendre transformation function of $\tilde f(y)$  defined in (\ref{5.7}).

In fact, consider the Legendre transformation of $\tilde{f}$
\begin{equation*}
x_i:=\tfrac{\partial \tilde{f}}{\partial y_i},\quad i=1,\ldots,n,\qquad u(x):=\sum^n_{i=1}y_{i}x_i-\tilde{f}(y),
\end{equation*}
then we obtain
\be \label{5.8}
x_1= -\tfrac{1}{4y_1}-\sum^n_{k=2}\tfrac{y^2_k}{2y_1^2},\qquad
x_j=\tfrac{ y_j}{ y_1},\quad 2\leq j\leq n.
\ee
and the function $u(x)$
\be\label{5.9}
u(x)=-\tfrac{1}{4}\bigg[\ln{\bigg(-x_1-\tfrac{1}{2}\sum^n_{k=2}x^2_k\bigg)}+\ln{4}+1 \bigg].
\ee
It is easy to see that $u(x)$ is Calabi-affine equivalent to $f(x)$ defined by (\ref{5.1}).

Similarly, we can prove that $\tilde{M}^n$ is complete with respect to the Calabi metric, and $\tilde f(y)$  satisfies the affine maximal type equation (\ref{1.2}) if and only if $a=-\frac{1}{n+1}$. If $a=-1$, then $\tilde f(y)$ satisfies the PDE (\ref{1.1}) with $L^\sharp=-4n(n+1)$. The proof of Theorem 1.3 is finished. \hfill $\Box$

\begin{rmk} The graph hypersurface of the function (\ref{4.4}) is also Euclidean complete and Calabi complete. However, it does not satisfy the PDE (\ref{1.1}).  The graph hypersurface of the function (\ref{4.1}) is neither Euclidean complete nor Calabi complete.
\end{rmk}

%%%%%%%%%%%%%%%%%%%%%%%%%%%%
%%%%%%%%%%%%%%%%%%%%%%%%%%%%
%%%%%%%%%%%%%%%%%%%%%%%%%%%
%%%%%%%%%%%%%%%%%%%%%%%%%%%%


\begin{thebibliography}{L3}

\bibitem{A} M. Abreu:  K\"{a}hler geometry of toric varieties and
extremal metrics. Int. J. Math. 9(1998), 641-651.

\bibitem{Ca} E. Calabi: Improper affine hyperspheres of convex type and a generalization of a theorem of J\"orgens. Michigan J. Math. 5(1958), 105-126.

\bibitem{CSL}L. Cao, L. Sheng, Z. Lian: A Bernstein property of certain types of fourth order partial differential equations. J. Math. Anal. Appl. 461(2018), 777-795.

\bibitem{CLS} B.H. Chen, A.-M. Li, L. Sheng: Affine techniques on extremal metrics on toric surfaces. Adv. Math.
340(2018), 459-527.

\bibitem{CH} X. Cheng, Z. Hu: Classification of locally strongly convex isotropic centroaffine hypersurfaces. Diff. Geom. Appl. 65(2019), 30-54.

\bibitem{CHM} X. Cheng, Z. Hu, M. Morus: Classification of the locally strongly convex centroaffine hypersurfaces with parallel cubic form. Results Math. 72(2017), 419-469.

\bibitem{D-2} S.K. Donaldson: Scalar curvature and stability
of toric varieties. J. Diff. Geom. 62(2002), 289-349.

\bibitem{D} S.Z. Du: Non-quadratic Euclidean complete affine maximal type hypersurfaces for
$\theta\in (0, (N-1)/N]$. J. Geom. Anal.  34(2024), no.8, Paper No. 229, 19 pp.

\bibitem{D2} S.Z. Du, X.Q. Fan: A Bernstein theorem for affine maximal-type hypersurfaces. C.R Acad. Sci. Paris. Ser. I, 357(2019), no. 1, pp. 66-73.

\bibitem{HLSV} Z. Hu, H. Li, U. Simon, L. Vrancken: On locally strongly convex affine hypersurfaces with parallel cubic form. Part I.  Diff. Geom. Appl. 27(2009), 188-205.

\bibitem{HLV} Z. Hu, H. Li, L. Vrancken: Locally strongly convex affine hypersurfaces with parallel cubic form. J. Diff. Geom. 87(2011), 239-307.

\bibitem {JL} F. Jia, A.-M. Li: Interior estimates
for solutions of a fourth order nonlinear partial differential
equation. Diff. Geom. Appl. 25(2007), 433-451.

\bibitem{LX} M.X. Lei, R.W. Xu: Classification of Calabi hypersurfaces in $\mathbb{R}^{n+1}$ with parallel Fubini-Pick tensor. J. Geom. Anal. 34(2024), no. 6, Paper No. 169, 30 pp.

\bibitem{LJ} A.-M. Li, F. Jia: A Bernstein property of some fourth order partial differential equations. Results Math. 56(2009), no.1-4, 109-139.

\bibitem{L-J-4} A.-M. Li, F. Jia:  A Bernstein property of affine
 maximal hypersurfaces. Ann. Global Anal. Geom.
 23(2003), 359-372.

\bibitem{LSZH} A.-M. Li, U. Simon, G.S. Zhao, Z.J. Hu: Global Affine Differential Geometry of Hypersurfaces. Second revised and extended edition, De Gruyter Expositions in Mathematics, 11, De Gruyter, Berlin, (2015).

\bibitem{LXSJ} A.-M. Li, R.W. Xu, U. Simon, F. Jia: Affine Bernstein Problems and Monge-Amp\`{e}re Equations. World Scientific Publishing Co. Pte. Ltd., Hackensack, NJ, 2010.

\bibitem{LXSJ1} A.-M. Li, R.W. Xu, U. Simon, F. Jia: Notes on Chern's affine Bernstein conjecture. Results Math.
60(2011), no. 1-4, 133-155.

\bibitem{LV}  H. Li, L. Vrancken: A basic inequality and new characterization of Whitney spheres in a complex
space form. Israel J. Math. 146(2005), 223-242.

\bibitem{L} Z. Lian: An invariant for affine maximal type equations.
J. Math. Anal. Appl. 543(2025), no. 2, part 1, Paper No. 128898, 12 pp.

\bibitem{LLX} X.Q. Liu, M.X. Lei, R.W. Xu: Classification of isotropic Calabi hypersurfaces. Diff. Geom. Appl.
85(2022), Paper No. 101954, 18 pp.

\bibitem{N1} S. N\"{o}lker: Isometric immersions of warped products, Diff. Geom. Appl., 6(1996), 1-30.

\bibitem{S} H. Shima: The geometry of Hessian structures. World Scientific Publishing Co. Pte. Ltd., Hackensack, NJ, 2007.

\bibitem{TW} N.S. Trudinger, X.-J. Wang:  The Bernstein problem
for affine maximal hypersurfaces. Invent. Math. 140(2000), 399-422.

\bibitem{TW1} N.S. Trudinger, X.J. Wang: Bernstein-J\"orgens theorem for a fourth order partial differential equation. J. Partial Differ. Equ. 15(2002), no.1, 78-88.

\bibitem{Xu-Cao} R.W. Xu, L.F. Cao: Bernstein properties for $\alpha$-complete hypersurfaces. J. Inequal. Appl. 2014, 2014:159, 9 pp.

\bibitem{XLL} R.W. Xu, A.-M. Li, X.X. Li: Euclidean complete $\alpha$ relative extremal hypersurfaces. J.
Sichuan University, 46(5)(2009), 1217-1223.

\bibitem{XL} R.W. Xu, X.X. Li: On the complete solutions to the Tchebychev-Affine-K\"ahler equation and its geometric significance. J. Math. Anal. Appl., 514(2022), no. 2, Paper No. 126351.

\bibitem{ZCH} Q. Zhao, X. Cheng, Z. Hu: On flat elliptic centroaffine Tchebychev hypersurfaces. J. Geom. Phys.
198(2024), Paper No. 105131, 15pp.

\end{thebibliography}
\end{document}